\documentclass[review,3p]{elsarticle}
\usepackage{mathrsfs,enumerate,amssymb,amsmath,bm,color,lineno}
\usepackage[all]{xy}
\usepackage{latexsym}
\usepackage{amsthm,a4wide,color}
\usepackage{amscd,verbatim}
\usepackage{graphicx}
\usepackage[colorlinks=true]{hyperref}

\modulolinenumbers[5]

\journal{Journal of \LaTeX\ Templates}

\newcommand{\pf}{\noindent\textbf{Proof.}\quad}
\newcommand{\epf}{\hspace{\stretch{1}}$\blacksquare$}

 \newtheorem{thm}{Theorem}[section]
 \newtheorem{cor}{Corollary}[section]
  
 \newtheorem{lem}{Lemma}[section]
 \newtheorem{prop}{Proposition}[section]
 \theoremstyle{definition}
 
 \newtheorem{ex}{Example}[section]
\theoremstyle{remark}
 
\begin{document}

\begin{frontmatter}

\title{Topological chain properties and shadowing property for dynamical systems on
uniform spaces\tnoteref{mytitlenote}}
\tnotetext[mytitlenote]{This work was supported by the National Natural Science Foundation of China (No. 11601449), Science and
Technology Innovation Team of Education Department of Sichuan for Dynamical System and its Applications (No. 18TD0013), and Youth
Science and Technology Innovation Team of Southwest Petroleum University for Nonlinear Systems (No. 2017CXTD02).}

\author[Ahmadi]{Seyyed Alireza Ahmadi}
\address[Ahmadi]{Department of Mathematics, University of Sistan and Baluchestan, Zahedan, Iran}
\ead{sa.ahmadi@math.usb.ac.ir, sa.ahmdi@gmail.com}

\author[Wu]{Xinxing Wu\corref{mycorrespondingauthor}}
\cortext[mycorrespondingauthor]{Corresponding author}
\address[Wu]{School of Sciences, Southwest Petroleum University, Chengdu, Sichuan 610500, P.R. China}
\ead{wuxinxing5201314@163.com}

\author[Chen]{Guanrong Chen}
\address[Chen]{Department of Electrical Engineering, City University of Hong Kong,
Hong Kong SAR, P.R. China}
\ead{gchen@ee.cityu.edu.hk}

\begin{abstract}
This paper discusses topological definitions of shadowing property, chain transitivity, totally chain transitivity,
and chain mixing property for dynamical systems on the uniform spaces and characterizes some topological chain property
for dynamical systems on compact uniform spaces. In particular, it is proved that a compact dynamical system is
topologically chain mixing if and only if it is totally topological chain transitive. Besides, some basic properties
for topological shadowing and non-wandering points on uniform spaces are obtained.
\end{abstract}

\end{frontmatter}


\section{Introduction}
The notion of shadowing plays an important role in the general qualitative theory of dynamical systems. In a dynamical system with the
shadowing property we can fined an orbit which remains close to a numerical solution for a long time. So dynamical systems with the shadowing
property behave nicely under computer simulations, in the sense that pseudo-orbits generated by the computer can be regarded as true orbits,
provided that the error of each step measured by a distance (a metric) on the space of calculation is sufficiently small.
Das et al.~\cite{DAS2013149} gave a purely topological definition for this notion, which is equivalent to the metric one in the case
that the phase space is a compact metric space. It is well known that many dynamical properties of dynamical systems with the shadowing property are equivalent.
In the other hand dynamical system with the shadowing property and the phase space inherit some properties to each others, for example the identity map on a
compact space has shadowing property if and only if the phase space is totally disconnected. A really convincing theory of topological dynamics exists only
with the assumption that the phase space $X$, in addition to being metric, is also compact. Most general results concerning chaotic properties, like positive
entropy, are obtained under this hypothesis  \cite{MR3539720,MOOTHATHU20112232,Wu2016,Wu2018TheCP,wu2019ijbc}. Sometimes one prefers to consider a
dynamical system on a non-metrizable topological space. On the other hand, many dynamical concepts in the topological theory of dynamical systems are
defined using the `distance' between sets and points.  However, for general topological spaces such distance -or size-related concepts cannot be
defined unless we have somewhat more structure than what the topology itself provides. This issue will be solved if we consider a completely regular,
and not necessarily metrizable, topological spaces which equipped with an structure, called uniformity, enabling us to control the distance between
points in these spaces.

Das et al.~\cite{DAS2013149} generalized the usual definitions in metric spaces of expansivity, shadowing, and chain recurrence for homeomorphisms
to topological spaces. We~\cite{WU2019145} introduced the topological concepts of weak uniformity, uniform rigidity, and multi-sensitivity and
obtained some equivalent characterizations of uniform rigidity. Then, we \cite{wu-ijbc2} proved that a point transitive dynamical system is either sensitive
or almost equicontinuous. Recently, we \cite{ahmadi3,ahmadi2} generalized concepts of entropy points, expansivity and shadowing property for dynamical systems
to uniform spaces and obtained a relation between topological shadowing property and positive uniform entropy. Shah et al. \cite{Shah2016}
showed that a dynamical system on a totally bounded uniform space which is topologically shadowing, mixing, and topologically expansive has the
topological specification property. Then, we \cite{ahmadi} proved that a dynamical system with ergodic shadowing is topologically chain transitive.
Good and Mac\'{\i}as \cite{good} obtained some equivalent characterizations and iteration invariance of various definitions of shadowing in the
compact uniformity sense generalizing the compact metric sense. For more results on shadowing
transitivity, and chain properties, one is referred to \cite{ahmadi4,ahmadi3,MOOTHATHU20112232,RICHESON2008251,WWL2018} and references therein.

This paper studies topological definitions of chain transitivity, totally chain transitivity, and chain mixing property
of dynamical systems with the topological shadowing property and obtains that topological chain recurrent, topological chain transitivity
and topological chain mixing properties are equivalent in a connected uniform space and use the topological shadowing property to characterize totally
disconnected uniform spaces.

\section{Basic definitions and preliminaries}\label{2}
A \textit{uniform space} is
a set with a uniform structure defined on it. A {\it uniform structure} $\mathscr{U}$ on a space $X$ is defined by the
specification of a system of subsets of the product $X\times X$ satisfying the following axioms:
\begin{itemize}
\item[U1)]
for any $E_1,E_2\in\mathscr{U}$, the intersection $E_1\cap
E_2$ is also contained in $\mathscr{U}$, and if $E_1\subset E_2$ and $E_1\in\mathscr{U}$, then
$E_2\in\mathscr{U}$;
\item[U2)] $\bigcap \mathscr{U} \supset \Delta_X  =
\{(x,x)~|~ x\in X\}$, i.e., every set $E\in\mathscr{U}$ contains the diagonal $\Delta_X$;
\item[U3)]
 if $E\in\mathscr{U}$, then $E^{\mathrm{T}} = \left\{(y,x)~|~(x,y)\in E\right\} \in\mathscr{U}$;
\item[U4)]
for any $E\in\mathscr{U}$, there exists $\hat{E}\in\mathscr{U}$ such that $\hat{E}\circ \hat{E}
\subset E$, where
$$
\hat{E}\circ \hat{E} =\left\{(x,y)~|~ \text{ there exists a } z\in
X \text{ such that }(x,z)\in \hat{E} \text{ and } (z,y)\in \hat{E}\right\}.
$$
\end{itemize}
The elements of $\mathscr{U}$ are called \textit{entourages} of the uniformity. A \textit{uniformity base}
or a \textit{basis of entourages} on a set $X$ is a
family of subsets of $X\times X$, such that the same three
conditions U2), U3) and U4) are satisfied, and which satisfies the conditions for a filter base (i.e. the
intersection of every pair of members of the family contains a member of the family). Given a
uniformity $\mathscr{U}$ on $X$ we can, of course, take the whole of $\mathscr{U}$ as a base, in
this sense. More usefully, we can take the symmetric entourages of $\mathscr{U}$ as a base, i.e.,
the entourages $E$ such that $E=E^\mathrm{T}$.
If $(X,\mathscr{U})$ is a uniform space, then the {\it uniform topology} on $X$ is
the topology in which a neighborhood base at a point $x\in X$ is formed by the family of sets $E[x]$,
where $E$ runs through the entourages of $X$, and $E[x]=\{y\in X~|~(x,y)\in E\}$ is called the
{\it cross section} of $E$ at $x$. By saying that a topological space is {\it uniformizable}
we mean, of course, that there exists a uniformity such that the associated uniform topology
is the given topology. It can be shown that a topological space
is uniformizable if and only if it is completely regular.

A mapping $f:X\longrightarrow Y$ from a uniform space $X$ into a uniform space $Y$
is {\it uniformly continuous} if the inverse image $(f\times f)^{-1}(E)$
is an entourage of $X$ for each entourage $E$ of $Y$. Throughout this paper,
assume that all uniform spaces are Hausdorff, i.e., $\bigcap \mathscr{U}=\Delta_{X}$, and
let $\mathbb{N}=\{1, 2, 3, \ldots\}$, $\mathbb{N}_0=\{0, 1, 2, \ldots\}$, and
$$
E^n:=E\circ E\circ\dots \circ E~~(n \textrm{ times}),
$$
for any $E\subset X\times X$.

A {\it dynamical system} (briefly, {\it system}) is a pair $(X, f)$, where $X$ is a uniform space
and $f:X\longrightarrow X$ is a uniformly continuous map.
For $x\in X$ and $U, V\subset X$, let
$$
N_f(x,U)=\{n\in\mathbb{N}_0~|~ f^n(x)\in U\} \mbox{ and } N_f(U,V)=\{n\in\mathbb{N}_0~|~U\cap f^{-n}(V)\neq\emptyset\}.
$$
An infinite subset $A\subset\mathbb{N}_0$ is {\it relatively dense} (or {\it syndetic}) if there exists $k>0$ such
that $\{n, n + 1\dots, n + k\}\cap A\neq\emptyset$ for all $n\in \mathbb{N}_0$ (this means that gaps are bounded);
and is {\it thick} if $A$ intersects every syndetic subset of $\mathbb{N}_0$.

Let $(X, f)$ be a dynamical system. A point $x\in X$ is
\begin{itemize}
\item[(1)] {\it recurrent} if $N_f(x, U)\setminus\{0\}\neq\emptyset$ for every neighborhood $U$ of $x$;
\item[(2)] {\it minimal} (or {\it almost periodic}) if $N_f(x, U)$ is syndetic for every neighborhood $U$ of $x$;
\item[(3)] {\it non-wandering} if $N_f(U, U)\setminus\{0\}\neq\emptyset$ for any neighborhood $U$ of $x$.
\end{itemize}
Denote the sets of all minimal, all recurrent, and all non-wandering points of $f$ by $M(f)$,
$R(f)$, and $\Omega(f)$, respectively. It is easy to see that $\Omega(f)$ is a closed invariant subset
of $X$, and $M(f)\subset R(f)\subset \Omega(f)$.

A dynamical system $(X, f)$ is
\begin{itemize}
\item[(1)] {\it non-wandering} if ${N_f(U, U)\setminus\{0\}}\neq\emptyset$ for every nonempty open subset $U$ of $X$;
\item[(2)] {\it transitive} if $N_f(U, V)\neq\emptyset$ for any pair of nonempty open subsets $U$, $V$ of $X$;
\item[(3)] {\it totally transitive} if $f^n$ is transitive for any $n\in\mathbb{N}$;
\item[(4)] {\it weakly mixing} if $f\times f$ is transitive.
\end{itemize}
Clearly, a dynamical system $(X, f)$ is non-wandering if and only if $\Omega(f)=X$.
Let $(X,f)$ be a dynamical system and let $D$ and $E$ be entourages of $X$. A {\it $(D,f)$-chain} of length $n$
is a sequence $\Gamma=\{x_i\}_{i=0}^{n}$  such that $(f(x_i),x_{i+1})\in D$ for $i=0, \ldots, n-1$. An infinite
$(D, f)$-chain is called a {\it $(D, f)$-pseudo-orbit}. If there exists no danger of confusion, we write $D$-chain
and $D$-pseudo-orbit instead of $(D, f)$-chain and $(D, f)$-pseudo-orbit, respectively. A $D$-pseudo-orbit
$\Gamma=\{x_i\}$ is {\it $E$-shadowed} by a point $y\in X$ if $(f^n(y),x_n)\in E$ for all $n \in \mathbb{N}_0$.
For $E\in\mathscr{U}$, use the symbol $\mathcal{O}_{E}(f,x,y)$ for the set of all $E$-chains $\{x_i\}_{i=0}^{n}$
of $f$ with $x_0=x$ and $x_n=y$ for some $n\in \mathbb{N}$. For any $x,y\in X$, we write $x\rightsquigarrow_E y$
if $\mathcal{O}_{E}(f,x,y)\neq\emptyset$ and write $x\rightsquigarrow y$ if $\mathcal{O}_{E}(f,x,y)\neq\emptyset$
for every $E\in\mathscr{U}$. We write $x\leftrightsquigarrow y$ if $x\rightsquigarrow y$ and $y\rightsquigarrow x$.
The set $\{x\in X~|~ x\leftrightsquigarrow x\}$ is called the \textit{chain recurrent set} of $f$ and is denoted by $CR(f)$.
A dynamical system $(X,f)$ has the {\it topological shadowing property} \cite{DAS2013149} if for every entourage $E$ of $X$, there exists
an entourage $D$ such that every $D$-pseudo-orbit is $E$-shadowed by some point $y$ in $X$. It is observed that for general topological
spaces, metric shadowing and topological shadowing are independent concepts \cite{DAS2013149}.
However, for compact metric spaces, metric shadowing and topological shadowing are equivalent.

A dynamical system $(X,f)$ is
\begin{itemize}
\item[(1)] {\it topologically chain transitive} if, for any entourage $E$ of $X$ and any two
points $x, y\in X$, $\mathcal{O}_{E}(f,x,y)\neq\emptyset$, or equivalently $x\rightsquigarrow_E y$;
\item[(2)] {\it topologically chain recurrent} if $CR(f)=X$;
\item[(2)] {\it totally topological chain transitive} if $f^n$ is topologically chain transitive for any $n\in\mathbb{N}$;
\item[(3)] {\it topologically chain mixing} if, for any two points $x,y\in X$ and any entourage $D$ of $X$,
there exists $N\in \mathbb{N}$ such that for any $n\geq N$, there exists a $D$-chain from $x$ to $y$ of length $n$.
\end{itemize}
A point $x\in X$ is an {\it equicontinuous point} for $f$, or say that $f$ is {\it equicontinuous} at $x$, if for every
entourage $E$ of $X$, there exists an entourage $D$ so that $(x,y)\in D$ implies $(f^n(x),f^n(y))\in E$ for all $n\in\mathbb{N}_0$.
A dynamical system $(X, f)$ is {\it equicontinuous} if, for every entourage $E$ of $X$, there exists an entourage $D$ such that $(f\times f)^n(D)\subset E$ for all $n\in\mathbb{N}_0$.
Note that in the case of compact metric spaces, for every entourage $E$ of $X$, there exists $\varepsilon>0$ such
that $d^{-1}([0,\varepsilon])\subset E$. Therefore, the above definitions coincide with the usual ones in compact
metric spaces. However, this does not hold if $X$ is not compact. For example, consider the entourage
$E=\{(x,y)\in\mathbb{R}^2~|~|x-y|<e^{-x^2}\}$ of $\mathbb{R}^2$. There exists no $\varepsilon>0$ such that
$d^{-1}([0,\varepsilon])\subset E$ (see \cite[Example~12]{DAS2013149}).

\section{Topological chain transitivity}
{This section is devoted to characterize some topological chain property for dynamical systems
on compact uniform spaces. In particular, it is proved that a compact dynamical system
is topologically chain mixing if and only if it totally topological chain transitive.}

\begin{lem}\label{Iterated-Lemma}
Let $(X, f)$ be a dynamical system on a uniform space $(X, \mathscr{U})$. Then,
$(X, f)$ has the topological shadowing property if and only if $(X, f^{n})$ has the 
topological shadowing property for
all $n\in \mathbb{N}$.
\end{lem}
\pf
From the definition of topological shadowing, this holds trivially.
\epf

\begin{lem}\label{L-E-Lemma}
Let $(X, f)$ be a topologically chain transitive system on a uniform space $(X, \mathscr{U})$
and $E\in \mathscr{U}$. Then, there exists $\iota_{E}\in \mathbb{N}$ such that, for any $y\in X$, $\iota_{E}$
is the greatest common divisor of the lengths of all $E$-chains from $y$ to itself.
\end{lem}
\pf
Fix $x \in X$ and let $\iota_{E}(x)$ be greatest common divisor of the lengths of $E$-chains from $x$ to itself.
For any $y \in X$ and any $E$-chain $\{y_{0} = y, y_{1}, \dots, y_{n} =y\}$ from $y$ to itself, since $(X, f)$ is topologically
chain transitive, there exists an $E$-chain $\{x_{0} = x, x_{1}, \dots, x_k=y, x_{k+1},\dots, x_{j}=x\}$ from $x$
to $y$ and back to $x$. Clearly, $j=m\iota_{E}(x)$ for some $m\in \mathbb{N}$.
Meanwhile, noting that the $E$-chain $\{x_{0} = x, x_{1}, \ldots, x_{k}=y, y_{1}, \ldots, y_{n} = y, x_{k+1}, \ldots, x_{j}=x\}$
from $x$ to itself has length $m\iota_{E}(x)+n$, one has $m\iota_{E}(x)+ n$ is necessarily a multiple of $\iota_{E}(x)$.
This implies that $\iota_{E}(x) $ divides $n$.
\epf

\begin{lem}\cite[Theorem 1.0.1]{Ramirez2005}\label{Prime}
Let $a_1, \ldots, a_n\in \mathbb{N}$. If $(a_1,\ldots, a_n)=1$, then there exists
$N\in \mathbb{N}$ such that any integer $s\geq N$ is representable as a non-negative integer
combination of $a_1, \ldots, a_n$.
\end{lem}

Let $f: (X, \mathscr{U})\longrightarrow (X, \mathscr{U})$ be a topological chain transitive map and $E\in \mathscr{U}$.
Define a relation `$\sim_{E}$' on $X$ by setting $x\sim_{E} y$ if and only if there exists an $E$-chain from
$x$ to $y$ of length a multiple of $\iota_{E}$.

\begin{itemize}
\item
Reflexive property: For any $x\in X$, it is clear that $x\sim_{E} x$;
\item
Transitive property: If $x\sim_{E} y$ and $y\sim_{E} z$, then by concatenating chains $x\sim_{E} z$;
\item
Symmetry property: If $x\sim_{E} y$, then there exists an $E$-chain $\Gamma_1$ from
$x$ to $y$ of length a multiple of $\iota_{E}$. Since $f$ is topologically chain transitive, there exists an $E$-chain $\Gamma_2$ from
$y$ to $x$ of length $m$. By concatenating $\Gamma_1$ and $\Gamma_2$, we obtain an $E$-chain from $x$ to $x$.
From Lemma~\ref{L-E-Lemma}, it follows that $m$ is a multiple of $\iota_{E}$. That is $y\sim_{E} x $.
\end{itemize}
Therefore, $\sim_{E}$ is an equivalent relation on $X$. If $x\sim_{E} y$, from the symmetry property of $\sim_{E}$,
it follows that there exists an $E$-chain $\Gamma'$ form $y$ to $x$ having length a multiple of $\iota_{E}$. Then,
for any $E$-chain $\Gamma$ from $x$ to $y$ of length $m$, concatenating $\Gamma$ and $\Gamma'$ obtains an $E$-chain
from $x$ to $x$. From Lemma~\ref{L-E-Lemma}, it follows that $m$ is a multiple of $\iota_{E}$. This mean that if
$x\sim_{E} y$, then every $E$-chain from $x$ to $y$ must have length a multiple of $\iota_{E}$.

For any $x\in X$, let $[x]_{E}$ be the equivalence class of $x$ for $\sim_{E}$. Since $(X, f)$ is topologically chain transitive,
there exists an $E$-chain $\Gamma$ from $f^{\iota_{E}}(x)$ to $x$ of length $m$. Clearly, $\Gamma_1:=\{x, f(x), \ldots, f^{\iota_{E}}(x)\}$
is an $E$-chain from $x$ to $f^{\iota_{E}}(x)$. Concatenating $\Gamma_1$ and $\Gamma$ obtains an $E$-chain from $x$ to $x$.
This implies that $m+\iota_{E}$ is a multiple of $\iota_{E}$, $m$ is as well. Thus, $f^{\iota_{E}}(x)\in [x]_{E}$, i.e.,
$f^{\iota_{E}}([x]_{E})=[x]_{E}$. For any $0\leq i<j<\iota_{E}$, it is clear that $f^{i}(x)\nsim_{E}f^{j}(x)$
(If $f^{i}(x)\sim_{E} f^{j}(x)$, then the length $j-i$ of chain $E$-chain $\{f^{i}(x), f^{i+1}(x), \ldots, f^{j}(x)\}$
must be a multiple of $\iota_{E}$, which is impossible).
Therefore, there exist $\iota_E$ equivalence classes for $\sim_{E}$, $f$ cycles among the classes periodically,
and every class is invariant under $f^{\iota_E}$. In fact, $X/\sim_{E}=\{[x]_{E}, [f(x)]_{E}, \ldots,
[f^{\iota_{E}-1}(x)]_{E}\}$ for any $x\in X$.

Similarly, we can define another equivalence relation `$\sim$' on $X$ by saying that
$x \sim y$ if and only if $x \sim_{E} y$ for all $E\in\mathscr{U}$. For any $x\in X$,
let $[x]$ denote the equivalence class of $x$ for $\sim$.

\begin{prop}\label{Open-Closed}
Let $(X, f)$ be a topologically chain transitive system on a uniform space $(X, \mathscr{U})$. Then,
for any $E \in \mathscr{U}$ and any $x\in X$, $[x]_{E}$
is both open and closed, and $[x]$ is closed.
\end{prop}
\pf
Choose $\hat{E}\in\mathscr{U}$ such that $\hat{E}^2\subset E$.
For any $y\in [x]_{E}$, since $f$ is topologically chain transitive, there exists an $\hat{E}$-chain
$\Gamma=\{x_0,\dots, x_{n}\}$ from $x$ to $y$. Noting that $\Gamma$ is also an $E$-chain, one has
$n$ is a multiple of $\iota_E$. For any $z\in\hat{E}[y]$, from $(f(x_{n-1}), y)\in \hat{E}$, it follows that
$(f(x_{n-1}), z)\in \hat{E}^{2}\subset E$. This means that $\{x_0, x_1, \ldots, x_{n-1}, z\}$ is an $E$-chain from
$x$ to $z$ having length a multiple of $\iota_{E}$, i.e., $\hat{E}[y]\subset [x]_{E}$. Thus, $[x]_{E}$ is open.

For any $y\in \overline{[x]_{E}}$, there exists $\hat{y}\in \hat{E}(y)\cap [x]_{E}$.
since $f$ is topologically chain transitive, there exists an $\hat{E}$-chain
$\hat{\Gamma}=\{\hat{x}_0,\dots, \hat{x}_{m}\}$ from $x$ to $\hat{y}$. Noting that $\hat{y}\in [x]_{E}$
and $\hat{\Gamma}$ is also an $E$-chain, one has $m$ is a multiple of $\iota_E$. From $(\hat{y}, y)\in \hat{E}$,
it follows that $(f(\hat{x}_{m-1}), y)\in \hat{E}^{2}\subset E$. This means that $\{\hat{x}_0, \hat{x}_1, \ldots, \hat{x}_{m-1}, y\}$
is an $E$-chain from $x$ to $y$ having length a multiple of $\iota_{E}$, i.e.,
$y\in [x]_{E}$. Thus, $\overline{[x]_{E}}=[x]_{E}$, i.e., $[x]_{E}$ is closed.
$[x]$ is also closed by applying $[x]=\bigcap_{E\in \mathscr{U}}[x]_{E}$.
\epf

\begin{lem}\label{l3}\cite[Proposition 8.16]{MR1687407}
Let $\mathscr{A}$ be an open covering of a compact uniform space $X$. Then there exists an
entourage $D$ such that $\mathfrak{C}(D)=\{D[x]~|~x\in X\}$ refines $\mathscr{A}$,
i.e., each of the uniform neighbourhoods $D[x]$
is contained in some member of $\mathscr{A}$.
\end{lem}

\begin{prop}\label{Le=1}
Let $(X, f)$ be a topologically chain transitive system on a compact uniform
space $(X, \mathscr{U})$ and $E\in\mathscr{U}$. If $f^{\iota_E}$ is topologically chain
transitive, then $\iota_E=1$.
\end{prop}
\pf
Suppose on the contrary that $\iota_{E}> 1$ and let $X=\bigcup_{i=1}^{\iota_{E}} A_{i} $, where $A_{i}$'s
are equivalence classes for $\sim_{E}$. By Proposition~\ref{Open-Closed}, the family
$\mathscr{A}=\{A_1,\dots, A_{\iota_{E}}\}$ is an open cover of $X$. This, together with Lemma~\ref{l3},
implies that there exists an entourage $V\subset E$ such that $\mathfrak{C}(V)=\{V[x]~|~x\in X\}$
refines $\mathscr{A}$. Fix $x\in A_i$ and $y\in A_j$ for $i\neq j$.
Since $f^{\iota_{E}}$ is topologically chain transitive, there exists a $(V, f^{\iota_{E}})$-chain
$\{x_{0} = x, x_{1}, \dots , x_{n} = y\}$ from $x$ to $y$. Since $(f^{\iota_{E}}(x_0),x_1)\in V$
and $f^{\iota_{E}}(x_0)\in A_i$, by construction of $V$ we have $x_1\in A_i$. Therefore,
by induction we conclude that $y=x_n\in A_i$, which is a contradiction.
\epf

\begin{lem}\label{Prime-PO}
Let $(X, f)$ be a topologically chain transitive system on a uniform space $(X, \mathscr{U})$ and $E\in \mathscr{U}$.
If $\iota_{E}=1$, then for any $x\in X$, there exist two $E$-chains from $x$ to itself with relatively prime lengths.
\end{lem}
\pf
Let $L_{E}:=\{n\in \mathbb{N}~|~ \text{there exists an } E\text{-chain from } x \text{ to } x \text{ of length } n\}
=\{n_1, n_2, \ldots, n_{k}, \ldots\}$ $(n_1<n_2<\cdots <n_k<\cdots)$. For any $k\in \mathbb{N}$, let $\Gamma_k$ be an $E$-chain
from $x$ to itself of length $n_{k}$. From proof of Lemma~\ref{L-E-Lemma}, it follows
that the greatest common divisor of $\{n_1, n_2, \ldots, n_{k}, \ldots\}$ is equal to $\min\{(n_1, n_2), (n_1, n_2, n_3), \ldots\}=1$.
This implies that there exist $n_1, n_2, \ldots, n_k\in L_{E}$ such that $(n_1, n_2, \ldots, n_k)=1$. Let $p=(n_1, n_2, \ldots, n_{k-1})$
and $p_{i}=\frac{n_i}{p}$ $(1\leq i\leq k-1)$. Then, $(p_1, p_2, \ldots, p_{k-1})=1$. By Lemma~\ref{Prime}, there exist
$\alpha_1, \alpha_2, \ldots, \alpha_{k-1}\in \mathbb{N}$ such that $(\alpha_1p_1+\alpha_2p_2+\cdots +\alpha_{k-1}p_{k-1}, n_{k})=1$.
By concatenating $\alpha_{1}$-$\Gamma_1$, $\alpha_2$-$\Gamma_2$, $\ldots$, and $\alpha_{k-1}$-$\Gamma_{k-1}$,
we obtain an $E$-chain $\Gamma$ from $x$ to $x$, i.e., $\Gamma=\underbrace{\Gamma_1\cdots \Gamma_1}\limits_{\alpha_1}
\underbrace{\Gamma_2\cdots \Gamma_2}\limits_{\alpha_2}\cdots \underbrace{\Gamma_{k-1}\cdots \Gamma_{k-1}}\limits_{\alpha_{k-1}}$.
Clearly, $\Gamma$ is an $E$-chain from $x$ to itself of length $q=(\alpha_1p_1+\alpha_2p_2+\cdots +\alpha_{k-1}p_{k-1})p$.
Clearly, $(q, n_{k})=1$ as $(p, n_{k})=1$. This means that $\Gamma$ and $\Gamma_{k}$ are two $E$-chains from $x$ to itself
with relatively prime lengths.
\epf

\begin{lem}\label{length<M}
Let $(X, f)$ be a dynamical system on a compact uniform space $(X, \mathscr{U})$. If $(X, f)$ is
topologically chain transitive, then for each $E\in\mathscr{U}$, there exists $M>0$ such that
for each $x,y\in X$, there exists an $E$-chain from $x$ to $y$ of
length less than or equal to $M$.
\end{lem}
\pf
Fix $E\in\mathscr{U}$ and choose an entourage $\hat{E}\in\mathscr{U}$ such that $\hat{E}^2\subset E$.
By compactness, there exist points $x_1,x_2,\dots,x_n\in X$ such that $X=\bigcup_{i=1}^{n}\hat{E}[x_i]$.
Since $(X, f)$ is topologically chain transitive, then for any $i,j\in\{1,2,\dots,n\}$, there exists an
$\hat{E}$-chain from $x_i$ to $x_j$ of length $m_{i,j}$. Put $M=\max \{m_{i,j}+1~|~i,j\in\{1,2,\dots,n\}\}$.
If $x,y\in X$, then there exist $k,l\in\{1,2,\dots,n\}$ such that $f(x)\in \hat{E}[x_k]$ and $y\in\hat{E}[x_l]$.
Choose an $E$-chain $\{x_k=z_0,z_1,\dots,z_{m_{k,l}}=x_l\}$ from $x_{k}$ to $x_{l}$ of length $m_{k, l}$.
Then, we obtain an $E$-chain $\{x,z_0,z_1,\dots,z_{m_{k,l}-1},y\}$ from $x$ to $y$ of length less than or equal to $M$.
\epf

\begin{thm}\label{1}
A dynamical system $(X, f)$ on a compact uniform space $(X, \mathscr{U})$
is totally topological chain transitive if and only if it is topologically chain mixing.
\end{thm}
\pf
By definition, if $f$ is topologically chain mixing, then it is totally topological chain transitive.

Fix any $E\in \mathscr{U}$. Assume that $f$ is totally topological chain transitive, applying
Proposition~\ref{Le=1} and Lemma~\ref{Prime-PO} yields that there exist two $E$-chains $\Gamma_1$
and $\Gamma_2$ from $x$ to itself with relatively prime lengths. By concatenating copies of
these loops and applying Lemma~\ref{Prime}, we can get an $E$-chain from $x$ to itself of length
$s$ for any $s> N$ or some $N\in \mathbb{N}$. By Lemma~\ref{length<M}, there exists some
$M\in \mathbb{N}$ such that between any two points in $X$, there exists an
$E$-chain of length less than or equal to $M$. By adding a loop at $x$ to a chain from $x$ to $y$,
we can get a $E$-chain from $x$ to $y$ of any length greater than $M + N$, implying that $f$ is
topologically chain mixing.
\epf

The proof of the following proposition is similar to Theorem \ref{1}.
\begin{prop}\label{p1}
Let $(X, f)$ be a dynamical system on a connected compact uniform space $(X, \mathscr{U})$.
If $(X, f)$ is topologically chain transitive, then it is topologically chain mixing.
\end{prop}
\pf
For any $E\in\mathscr{U}$, we show that $\iota_E=1$. If $\iota_E >1$, then $X=\bigcup_{i=1}^{\iota_{E}} A_{i}$,
where $A_{i}$'s are equivalence classes with respect to $\sim_E$, which are closed and open subsets of $X$ by
Proposition~\ref{Open-Closed}, which contradicts the connectivity of $X$. An argument similar to the proof of
Theorem~\ref{1} shows that $f$ is topologically chain mixing.
\epf

\begin{cor}\label{c1}
Let $(X, f)$ be a dynamical system on a connected  compact uniform space $(X, \mathscr{U})$.
Then, the following statements are equivalent:
\begin{enumerate}[(1)]
\item $f$ is topologically chain recurrent;
\item $f$ is topologically chain transitive;
\item $f$ is totally topological chain transitive;
\item $f$ is topologically chain mixing.
\end{enumerate}
\end{cor}
\pf
Clearly $(4) \Longrightarrow (3) \Longrightarrow (2) \Longrightarrow (1)$. By Proposition~\ref{p1},
$(2)\Longrightarrow (4)$. It suffices to prove that $(1)\Longrightarrow (2)$.

Suppose that $f$ is topologically chain recurrent and $E\in\mathscr{U}$. We say that $x$ and $y$
are $E$-chain equivalent if $x\rightsquigarrow_E y$ and $y\rightsquigarrow_E x$. Since $f$ is topologically
chain recurrent, this is an equivalence relation. Let $x$ and $y$ be $E$-chain equivalent. Choose
an entourage $D\in \mathscr{U}$ such that $D^2\subset E$ and $(f\times f)(D)\subset E$.
Now we show that $x$ is $E$-chain equivalent to every $z\in D[y]$. Let $\{x_0=x,x_1,\dots,x_n=y\}$
be an $E$-chain from $x$ to $y$ and $\{z_0=y,z_1,\dots ,z_m=y\}$ be an $D$-chain from $y$ to itself.
Then $\{x_0=x,x_1,\dots,x_n=y=z_0,z_1,\dots,z_m=z\}$ is an $E$-chain from $x$ to $z$. Similarly,
let $\{y=y_0,y_1,\dots,y_p=x\}$ be an $E$-chain from $y$ to $x$. Then,
$\{z,z_1,\dots, z_m=y=y_0,y_1,\dots,y_n=x\}$ is an $E$-chain from $z$ to $x$. Thus,
$x$ is $E$-chain equivalent to $z$, implying that the equivalence class $[x]$ is open and therefore closed.
This, together with connectivity of $X$, implies that $X=[x]$. Therefore, $f$ is topologically chain transitive.
\epf

The next example shows that Corollary \ref{c1} does not hold for the non-connected case.
\begin{ex}\label{ex1}\cite{ahmadi}
Let $\mathbb{P}=\mathbb{R}\setminus \mathbb{Q}$. Suppose that $\textbf{a}=\{a_i\}_{i\in\mathbb{Z}}\subset\mathbb{P}$
is an increasing  bi-sequence for which there exists a positive integer $k$ such that $a_i+1=a_{i+k}$ for all
$i\in\mathbb{Z}$. Put
$$
U_{\mathbf{a}}=\Cup_{i\in\mathbb{Z}}[\{(a_i,a_i)\}\cup (a_i,a_{i+1})\times(a_i,a_{i+1})]
$$
\begin{figure}[ht]
\centerline{\includegraphics[scale=0.12]{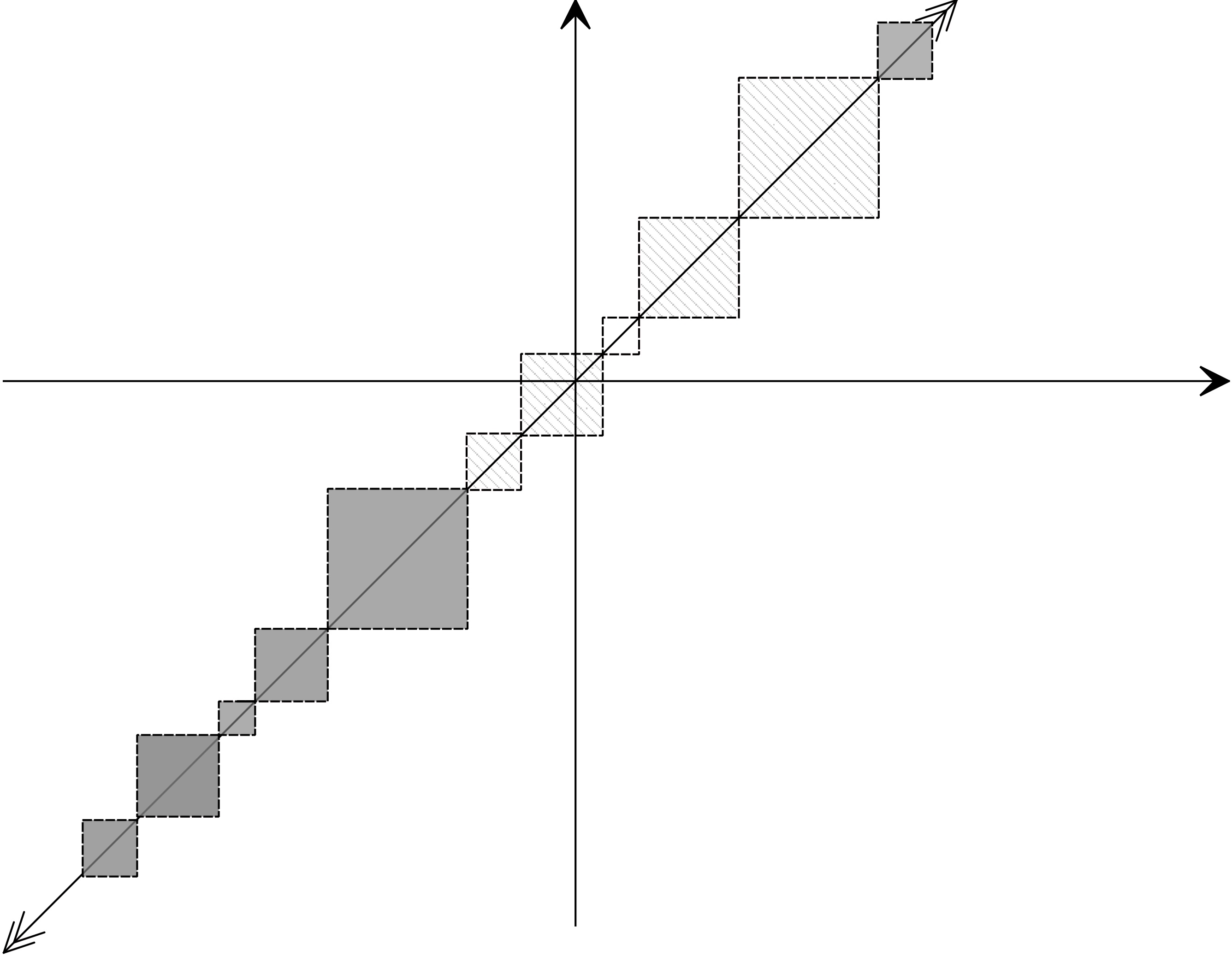}}
\caption{ Filter base for Michael line uniformity}
\end{figure}
(See Figure~1). Then the family $\mathcal{B}=\{U_{\textbf{a}}\}$ is a filter base and the uniformity generated
by this filter base generate a topology in which any point of $\mathbb{P}$ is isolated.
Let $\mathbb{S}^1$ be the unit circle. Consider the uniformity $\mathscr{U}$ on $\mathbb{S}^1$ -- just taking the
projection modulo $1$ -- and define $f:\mathbb{S}^1\longrightarrow \mathbb{S}^1$ by
$f(x)=x+\alpha$, $\alpha\in\mathbb{P}$. Then, it can be verified that $(\mathbb{S}^1,f)$ is totally
topological chain transitive, topologically chain mixing and topologically chain recurrent,
which does not have topological shadowing and mixing property.
\end{ex}

\section{Topological shadowing property}
This section obtains some basic properties for non-wandering points set $\Omega(f)$. In particular,
it is proved that a compact dynamical system $(X, f)$ share the same topological shadowing
property with its restricted system $(\Omega(f), f|_{\Omega(f)})$.

\begin{prop}\label{Infinite-Omega}
Let $(X, f)$ be a dynamical system on a uniform space $(X, \mathscr{U})$ and $\Omega(f)\neq \emptyset$.
Then, for any encourage $U\in \mathscr{U}$, $N_{f}(U[x], U[x])$ is an infinite set.
\end{prop}

\pf
If $x$ is a periodic point, then clearly $N_{f}(U[x], U[x])$ is infinite. Suppose $x$ is not a periodic point.
Then, for any fixed $N\in \mathbb{N}$, $f^{i}(x)\neq x$ for all $1\leq i\leq N$. Since $(X, \mathscr{U})$ is
Hausdorff, there exists an encourage $\hat{U}\subset U$ such that $\hat{U}[f^{i}(x)]\cap \hat{U}[x]=\emptyset$
for all $1\leq i\leq N$. This, together with the continuity of $f$, yields that there exists an encourage
$D\subset \hat{U}$ such that, for any $1\leq i\leq N$,
$$
f^{i}(D[x])\subset \hat{U}[f^{i}(x)],
$$
implying that
$$
f^{i}(D[x])\cap D[x]\subset \hat{U}[f^{i}(x)]\cap \hat{U}[x]=\emptyset.
$$
Then, there exists $n>N$ such that $f^{n}(D[x])\cap D[x]\neq \emptyset$ by $x\in \Omega(f)$.
\epf

\begin{prop}
Let $(X, f)$ be a dynamical system on a uniform space $(X, \mathscr{U})$ and $\Omega(f)\neq \emptyset$.
If $(X, f)$ has the topological shadowing property and $x\in\Omega(f)$, then for any entourage $U\in \mathscr{U}$,
there exists a natural number $k$ and a point $w\in U[x]$
such that $f^{nk}(w)\in U[x]$ for all $n\in\mathbb{N}$.
\end{prop}
\pf
Choose an encourage $E\in\mathscr{U}$ such that $E^2\subset U$ and take an entourage
$D\subset E$ such that every $D$-pseudo-orbit can be $E$-shadowed by some point in $X$
by the topological shadowing property. Pick an encourage $\hat{D}\in \mathscr{U}$ with $\hat{D}^2\subset D$.
Since $x\in \Omega(f)$, there exists $z\in \mathrm{Int}(\hat{D}[x])\cap f^{-k}(\mathrm{Int}(\hat{D}[x]))$
for some $k\in\mathbb{N}$. If we write $\Gamma=\{z,f(z),\dots, f^{k-1}(z)\}$, then $\Gamma\Gamma\Gamma\cdots$
is an infinite $D$-pseudo-orbit of $f$ as $(f^{k}(z),z)\in \hat{D}^2\subset D$. Let $w\in X$ be a
point which $E$-shadows this pseudo-orbit. Then, $(z,f^{nk}(w))\in E$ for all $n\in \mathbb{N}_0$.
This, together with $(x, z)\in \hat{D}$, implies that $(x,f^{nk}(w))\in \hat{D}\circ E\subset E^2\subset U$, i.e.,
$f^{nk}(w)\in U[x]$.
\epf

\begin{cor}\label{c2}
Let $(X, f)$ be a dynamical system on a uniform space $(X, \mathscr{U})$ and $\Omega(f)\neq \emptyset$.
If $(X, f)$ has the topological shadowing property and $x\in\Omega(f)$, then for any open set $U$ containing
$x$, $k\mathbb{N}_0\subset N_f(U,U)$ for some $k\in \mathbb{N}$.
\end{cor}

\begin{prop}\label{lem2}
Let $(X, f)$ be a dynamical system on a uniform space $(X, \mathscr{U})$.
If $(X, f)$ has the topological shadowing property
and totally transitive, then it is weakly mixing.
\end{prop}
\pf
Since $f$ is totally transitive, it is clear that $\Omega(f)\neq \emptyset$.
Given any two nonempty open subsets $U$ and $V$ of $X$, from Corollary \ref{c2},
it follows that there exists $k\in\mathbb{N}$ such that $k\mathbb{N}\subset N_f(U,U)$.
Since $f^k$ is transitive, $k\mathbb{N}_0\cap N_f(U,V)\neq\emptyset$. Therefore,
$N_f(U,U)\cap N_f(U,V)\neq\emptyset$,
implying that $f$ is weakly mixing.
\epf

The {\it $\omega$-limit set} $\omega(x, f)$ of $x$ consists of all
the limit points of $\{f^{n}(x)~|~n\in \mathbb{N}_0\}$, i.e.,
$$
\omega(x, f)=\{y\in X~|~\exists n_{k}\nearrow +\infty \text{ with } f^{n_k}(x)\rightarrow y\}.
$$
Let $\omega(f)=\bigcup_{x\in X}\omega(x, f)$. Clearly, $\omega(f)\subset \Omega(f)$ and $\omega(f)\neq\emptyset$
for every compact dynamical system $(X, f)$.

\begin{prop}\label{Minimal}
Let $(X, f)$ be a dynamical system on a compact uniform space $(X, \mathscr{U})$.
If $(X, f)$ has the topological shadowing property and $x\in\Omega(f)$, then for
any entourage $U\in \mathscr{U}$, there exists a natural number $k$ and a point $y\in U[x]\cap M(f)$
such that $f^{nk}(y)\in U[x]$ for all $n\in\mathbb{N}$.
\end{prop}
\pf
Choose an encourage $E\in\mathscr{U}$ such that $E^3\subset U$ and take an entourage
$D\subset E$ such that every $D$-pseudo-orbit can be $E$-shadowed
by some point in $X$. Pick an encourage $\hat{D}\in \mathscr{U}$
with $\hat{D}^2\subset D$. Since $x\in \Omega(f)$, there exists
$z\in \mathrm{Int}(\hat{D}[x])\cap f^{-k}(\mathrm{Int}(\hat{D}[x]))$
for some $k\in\mathbb{N}$. If we write $\Gamma=\{z,f(z),\dots, f^{k-1}(z)\}$, then
$\Gamma\Gamma\Gamma\cdots$ is an infinite $D$-pseudo-orbit of $f$ as $(f^{k}(z),z)\in \hat{D}^2\subset D$.
Let $w\in X$ be a point which $E$-shadows this pseudo-orbit. Then, $(z,f^{nk}(w))\in E$ for all $n\in \mathbb{N}_0$.
This, together with $(x, z)\in \hat{D}$, implies that $(x,f^{nk}(w))\in \hat{D}\circ E\subset E^2$ for all $n\in \mathbb{N}_0$.
Let $y\in\overline{\{f^{nk}(w)~|~n\in \mathbb{N}_0\}}$ be a minimal point for $f^{k}$ (existence is by Zorn's lemma and the
compactness of $X$). Clearly, $y\in M(f)$. For any $m\in \mathbb{N}_0$, there exists $n_{m}\in \mathbb{N}_0$ such that
$(f^{n_mk}(w), f^{mk}(y))\in E$, implying that $(x, f^{mk}(y))\in E^2 \circ E=E^3\subset U$, i.e., $f^{mk}(y)\in U[x]$.
\epf

The following is an immediate corollary of Proposition~\ref{Minimal} and \cite[lemma 2.8]{akin2001}
which states that if the minimal points for $f$ and $g$ are dense in $X$ and $Y$ respectively,
then the minimal points of $f\times g$ are dense in $X\times Y$.
\begin{cor}\label{Omega=M}
Let $(X, f)$ be a dynamical system on a compact uniform space $(X, \mathscr{U})$.
If $(X, f)$ has the topological shadowing property,
then $\Omega(f)^n=\overline{M(f^{(n)})}$ for any $n\in \mathbb{N}$.
\end{cor}

A uniform $(X, \mathscr{U})$ is
\begin{enumerate}[(1)]
\item {\it connected} if $X$ contains no open and closed subset, other than the empty
set and the full set;
\item {\it totally disconnected} if all the connected subspaces of $X$ are one-point sets.
\end{enumerate}
We can define an equivalence relation `$\sim$' on $X$ by setting $x\sim y$ if there
exists a connected subspace of $X$ containing both $x$ and $y$. The equivalence classes are
called the {\it connected components} of $X$, or simply {\it components}.
Clearly, $X$ is totally disconnected if and only if each connected component of $X$ is a singleton.
Recall that a locally compact Hausdorff topological space $X$ is {\it totally disconnected}
if and only if it has a basis of topology consisting of compact open sets.

Suppose that $\mathscr{B}$ is a basis of entourages for a uniform space $X$. An injection $f$ of $X$ into itself is an {\it isobasism with respect to $\mathscr{B}$} if,
for every entourage $V$ of $\mathscr{B}$, $(x, y)\in V$ is equivalent to $(f(x), f(y)) \in V$.
When $f$ maps $X$ on itself, this condition is
equivalent to $(f\times f) (V)= V$.

\begin{lem}\label{isom}\cite[Theorem~8]{rhodes_1956}
If $G$ is a uniformly equicontinuous group of maps of a uniform space $X$ on itself,
then there exists a basis of entourages $\mathscr{B}$ on $X$ such that every map
$f$ of $G$ is a $\mathscr{B}$-isobasism.
\end{lem}

\begin{thm}
Let $(X,f)$ be a surjective equicontinuous dynamical system on a compact uniform space $(X, \mathscr{U})$.
Then, $(X, f)$ has the topological shadowing property if and only if $X$ is totally disconnected.
\end{thm}
\pf
By Lemma \ref{isom},
there exists a basic neighborhood $\mathscr{B}$ on $X$ such that $f$ is $\mathscr{B}$-isobasism.

First assume that $f$ has the topological shadowing property.
Suppose on the contrary that $X$ is not totally disconnected. Then, $X$ has a non-degenerate component $C$.
Let $x$ and $y$ be two distinct points of $C$. Since $X$ is Hausdorff, there exists an encourage $E\in\mathscr{B}$
such that $(x,y)\notin E$. Choose $\hat{E}\in\mathscr{B}$ such that $\hat{E}^2\subset E$.
By the topological shadowing property, there exists $D\in\mathscr{B}$ such that every $D$-pseudo-orbit
can be $\hat{E}$-shadowed by some point in $X$. Since $C$ is connected, there exists a sequence
$\{x=x_0,x_1\dots,x_n=y\}$ in $C$ such that $(x_i,x_{i+1})\in D$ for all $0\leq i \leq n-1$.
Let $y_i=f^i(x_i)$ ($0\leq i \leq n-1$). Then, for any $0\leq i \leq n-1$,
we have $(f^{i+1}(x_{i+1}),f^{i+1}(x_i))\in D$ by applying Lemma~\ref{isom}, implying that $(y_{i+1},f(y_i))\in D$.
That means that $\{y_0, y_1, \dots, y_n\}$ is a $D$-pseudo-orbit and there exists a point $z\in X$
such that $(y_i,f^i(z))\in\hat{E}$ for all $0\leq i \leq n-1$. Since $f$ is $\mathscr{B}$-isobasism
and $(f^n(z),y_n)=(f^n(z),f^n(x_n))\in\hat{E}\in \mathscr{B}$, then $(z, x_{n})\in \hat{E}$. This, together with
$(z,x_0)=(z,y_0)\in \hat{E}$, implies that
it follows that $(x,y)=(x_0,x_n)\in\hat{E}\circ\hat{E}\subset E$, which is a contradiction.

Conversely assume that $X$ is totally disconnected. Given any entourage $E\in\mathscr{U}$ and choose
$\hat{E}\in\mathscr{U}$ such that $\hat{E}^2\subset E$. Since $X$ is
totally disconnected, then, for any $z\in X$, there exists an open and compact
subset $W_z$ of $X$ such that $z\in W_z\subset \mathrm{Int}(\hat{E}[z])\subset \hat{E}[z]$.
By the compactness of $X$, there exist $z_1, \dots, z_m\in X$ such that $X=\bigcup_{j=1}^{m}W_{z_j}$.
Applying Lemma \ref{l3} yields that there exists an encourage $D\in\mathscr{B}$ with $D\subset E$
which refines the open cover $\{W_{z_1}, W_{z_2}, \ldots, W_{z_m}\}$. Let
$\{x_0,x_1,\dots\}$ be a $D$-pseudo-orbit of $f$. For any fixed $n\in \mathbb{N}_0$,
from $(f(x_n),x_{n+1})\in D$, it follows that there exists $1\leq i\leq m$ such that
$\{f(x_n),x_{n+1}\}\subset W_{z_i}$. This, together with $(f^2(x_{n-1}),f(x_n))\in D$
(by the definition of isobasism), implies that $f^2(x_{n-1})\in W_{z_i}$.
Going backward inductively, we have $f^{n+1}(x_0)\in W_{z_i}$. Thus,
$\{f^{n+1}(x_0),x_{n+1}\}\subset W_{z_i}\subset\hat{E}[z_1]$,
implying that $(f^{n+1}(x_0),x_{n+1})\in E$. This means that
every $D$-pseudo-orbit of $f$ can be $E$-shadowed by the starting point.
Therefore, $(X, f)$ has the topological shadowing property.
\epf

Let $(X,f)$ be a compact dynamical system on a uniform space $(X, \mathscr{U})$ and $D\in \mathscr{U}$.
Define a relation $\backsimeq_{D}$ on $\Omega(f)$ by $x\backsimeq_D y$ if and only if there exist
$D$-chains  from $x$ to $y$ and from $y$ to $x$ in $\Omega(f)$. From the definition of non-wandering points,
it can be verified that for any $x\in \Omega(f)$ and any $E\in \mathscr{U}$, there exists an $E$-chain
from $x$ to itself. This implies that $\backsimeq_{D}$ is an equivalence relation on $\Omega(f)$.

For any $x\in \Omega(f)$, let $[x]_{D}$ be the equivalence class of $x$ for $\backsimeq_{D}$ in $\Omega(f)$. For any $y\in [x]_{D}$,
there exist $D$-chains $\Gamma_1$ and $\Gamma_2$ in $\Omega(f)$ such that $\Gamma_1=\{x_0=x, x_1, \ldots, x_{n-1}, x_{n}=y\}$
and $\Gamma_2=\{y_0=y, y_1, \ldots, y_{m-1}, y_{m}=x\}$. From $y\in D[f(x_{n-1})]\cap f^{-1}(D[y_1])$,
it follows that there exists an encourage $\hat{D}\subset D$ such that $\hat{D}[y]\subset D[f(x_{n-1})]\cap f^{-1}(D[y_1])$.
Clearly, for any $z\in \hat{D}[y]\cap \Omega(f)$, $\hat{\Gamma}_1=\{x, x_1, \ldots, x_{n-1}, z\}$ and
$\hat{\Gamma}_2=\{z, y_1, \ldots, y_{m-1}, x\}$ are $D$-chains in $\Omega(f)$,
implying that $\hat{D}[y]\cap \Omega(f)\subset [x]_{D}$. Thus, $[x]_{D}$ is open in $\Omega(f)$.

Take an encourage $\hat{D}\subset D$ such that $(f\times f)^{2}(\hat{D})\subset D$. From Proposition~\ref{Infinite-Omega},
it follows that there exists $n>2$ such that $f^{n}(\hat{D}[x])\cap \hat{D}[x]\neq \emptyset$. Then, there exists $z\in \hat{D}[x]$
such that $f^{n}(z)\in \hat{D}[x]$. Clearly, $(f^{2}(z), f^{2}(x))\in D$ by the choice of $\hat{D}$. This implies that
$\{f(x), f^{2}(z), \ldots, f^{n-1}(z), x\}$ is a $D$-chain from $f(x)$ to $x$. Clearly, $\{x, f(x)\}$ is a $D$-chain from $x$
to $f(x)$. Thus, $x\backsimeq_{D}f(x)$, i.e., $f(x)\in [x]_{D}$. Similarly, it can be verified that $f^{n}(x)\in [x]_{D}$
for all $n\geq 2$. This implies that, for any $y\in [x]_{D}$, $f(y)\in [y]_{D}=[x]_{D}$, i.e, $f([x]_{D})\subset [x]_{D}$.

\begin{thm}
Let $(X,f)$ be a dynamical system on a compact uniform space $(X, \mathscr{U})$. If $(X, f)$ has the
topological shadowing property, then $(\Omega(f), f|_{\Omega(f)})$ has the topological shadowing property.
\end{thm}
\pf
For any fixed $E\in \mathscr{U}$, choose $\hat{E}\in\mathscr{U}$ such that $\hat{E}^3\subset E$. From the topological shadowing
property of $(X, f)$, it follows that there exists an encourage $D\subset \hat{E}$ such that every $D$-pseudo-orbit
can be $\hat{E}$-shadowed by some point in $X$. From the above discussion and compactness of $\Omega(f)$,
there exist $x_1, x_2, \ldots, x_{m}\in \Omega(f)$ such that $\Omega(f)=\bigcup_{i=1}^{m}[x_i]_{D}$
and $[x_1]_{D}, \ldots, [x_m]_{D}$ are mutually disjoint. Let $\hat{D}\subset D$ be the entourage from Lemma~\ref{l3}
which refines the open cover $\{[x_1]_{D}, \dots, [x_m]_{D}\}$.

{\bf Claim.} Every $\hat{D}$-pseudo-orbit in $\Omega(f)$
can be $E$-shadowed by some point in $\Omega(f)$.

Let $\{z_{s}\}_{s=0}^{+\infty}$ be a $\hat{D}$-pseudo-orbit in $\Omega(f)$. Then, there exists
$1\leq i\leq m$ such that $z_{0}\in [x_{i}]_{D}$. This, together with $f(z_{0})\in [x_i]_{D}$
and $(f(z_0), z_1)\in \hat{D}$, implies that $z_{1}\in [x_{i}]_{D}$ as $\hat{D}$
refines the open cover $\{[x_1]_{D}, \dots, [x_m]_{D}\}$. Similarly, it can be verified that
$z_{n}\in [x_{i}]_{D}$ for all $n\geq 2$.

Fix $s\in \mathbb{N}$. Clearly, $\{z_0, \dots, z_s\}$ is a $\hat{D}$-chain in $\Omega(f)\cap [x_i]_{D}$.
Since $[x_i]_{D}$ is an equivalence class, there exists a $D$-chain $\{y_0, \dots, y_t\}$ from $z_s$
to $z_0$. Let $\Gamma=\{z_0, \dots, z_{s-1}, y_0, \dots, y_{t-1}\}$. Then, $\hat{\Gamma}=\Gamma\Gamma\Gamma\cdots$
is a $D$-pseudo-orbit in $\Omega(f)$. 
Thus, there exists ${\bm z}_{s}\in X$ which $\hat{E}$-shadows $\hat{\Gamma}$. Put $k=s+t$,
then $(f^{kn+j}({\bm z}_{s}),z_j)\in \hat{E}$ for all $n\geq 0$ and $0\leq j<k$. Choose
$\hat{{\bm z}}_{s}\in \omega({\bm z}_{s}, f^{k})\subset \Omega (f)$
($\omega({\bm z}_{s}, f^{k})\neq \emptyset$ by the compactness of $X$). Since $f$ is uniformly continuous,
there exists an encourage $U\subset \hat{E}$ such that $(f\times f)^{j}(U)\subset \hat{E}$ holds for all
$0\leq j\leq k$. Take $n_{0}\in \mathbb{N}$ such that $(f^{kn_0}({\bm z}_{s}), \hat{{\bm z}}_{s})\in U$.
Then, for any $0\leq j\leq s$, $(f^{kn_0+j}({\bm z}_{s}), f^{j}(\hat{{\bm z}}_{s}))\in \hat{E}$ and
$(f^{kn_0+j}({\bm z}_{s}), z_{j})\in \hat{E}$, implying that
\begin{equation}\label{e-1}
(f^{j}(\hat{{\bm z}}_{s}), z_j)\in \hat{E}^{2}.
\end{equation}
Choose a limit point $\hat{{\bm z}}$ of the sequence
$\{\hat{{\bm z}}_{1}, \hat{{\bm z}}_{2}, \ldots, \hat{{\bm z}}_{s}, \ldots\}$.
Clearly, $\hat{{\bm z}}\in \Omega(f)$. Similarly, by the uniform continuity of $f$ and (\ref{e-1}),
it can be verified that, for any $s\in \mathbb{N}_0$,
$(f^{s}(\hat{{\bm z}}), z_{s})\in E$.
\epf

\section*{References}

\bibliography{mybibfile}

\end{document}